# Digital-Discrete Surface Reconstruction: A true universal and nonlinear method

Li Chen
Department of Computer Science and Information Technology
University of the District of Columbia
Email: lchen@udc.edu
March 9, 2010

**Abstract**
The most common problem in data reconstruction is to fit a function based on the observations of some sample (guiding) points. This paper provides a methodological point of view of digital-discrete surface reconstruction. We explain our method along with why it is a truly universal and nonlinear method unlike most popular methods, which are linear and restricted. This paper focuses on what the surface reconstruction problem is and why the digital-discrete method is important, necessary, and how it can be accomplished.

## 1. Surface Reconstruction and Digital-Discrete Reconstruction

The most common problem in data reconstruction is fitting a function based on the observations of some sample (guiding) points. If the function is defined as a polynomial or a constraint such as a partial differential equation (PDE), we can try to find the coefficients of the polynomial or find a solution to the equation. The most popular methods for answering these questions include the Splines method for polynomials and the finite elements method for PDE. The finite difference method is usually used when the boundary is provided and uses a system of linear equations. It is not usually used in general cases since guiding points are often unavailable for regular grid points.

Formally, for randomly arranged sample points on a plane or manifold, we can assume a sample data set S={ (p, v(p)) | p is a point on the manifold, v(p) is a real number} where S is finite. There are two ways to solve the problem: (1) Find the analytical solution, i.e. we find the exact solution if it exists and can be expressed; (2) Find the numerical solution where we can usually only get an approximation of the solution. In most cases, applications are not solvable analytically. As a result, the primary focus is on the numerical solution.

The classical methods for this problem are divided into two categories: mesh-dependent and mesh-independent [3][1]. The mesh-dependent method requires convex-cell decomposition, usually a triangulation, on the domain based on sample points, e.g. the finite element method. The mesh-independent method requires a dense sample point set and has become a very popular topic in recent years especially in computer graphics [19][2][16]. A common method used today is called the moving least square method, which uses a circular weighted neighborhood instead of triangles [24][21][17]. The advantage of these methods is that they can do localized fittings. On the contrary, the finite difference method can only do global fittings, meaning that every sample point must be considered for each undetermined new data point. Global fitting using the finite element method or moving least square method is too costly in practice. A definite disadvantage of using both of these methods is that they would not be able to deal with a very small set of sample points.

In [8], a systematic digital-discrete method was designed for obtaining continuous functions with smoothness of a certain order ($C^n$) from sample data. This method has the advantage of being able to use neighbor sample points without decomposition. Also, it can consider all sample points or just do a local fitting. The new method is nonlinear and based on discrete mathematics especially graph theory. It overcomes the problem of the traditional discrete decomposition method, which relies on a specific cell decomposition of the domain. The new method also limits the problem that a circular neighborhood based method has for a closed manifold, e.g. it may not work correctly if "none-true" neighbors are selected when the sample data set is not dense enough for some unfitted points.

Our method is based on gradually varied functions [4][5][6][7]; see Fig.1. The new algorithm obtains the best solution to the fitting. We have added a component of the classical finite difference method for smooth reconstruction. The new digital-discrete method has considerable advantages in a large number of real data applications. This method can potentially be used to obtain smooth functions such as polynomials through its derivatives $f^{(k)}$ and the solutions for partial differential equations including harmonic and other important equations.

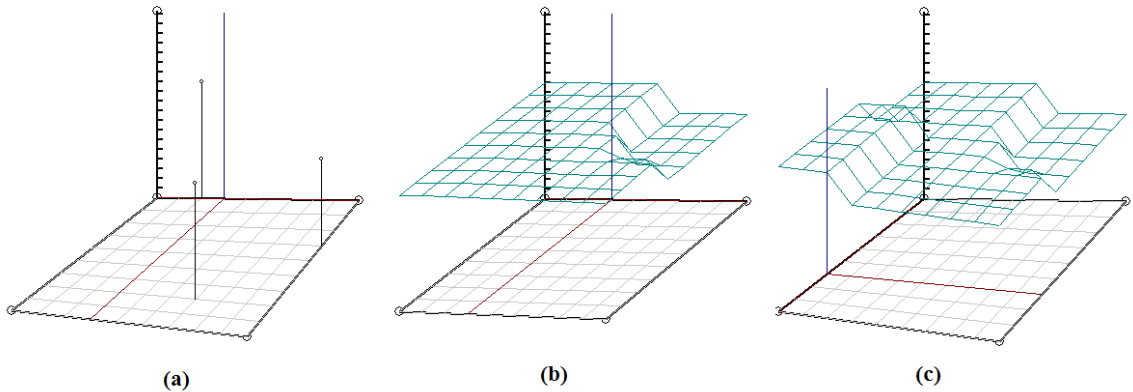

(a)    (b)    (c)

Fig 1. Gradually varied interpolations: (a) Three guiding points, (b) Fitting the points, (c) One more guiding point added on the left side.

This paper provides a thorough explanation of why the digital-discrete method is necessary. It can be considered as a methodology paper comparing the applications of this method with the method described in [10]. Our conclusion is that this method is a truly universal method meaning that it works for both local and global fittings. Using as many sample points as we would like, it can still fit the local region regardless of the other unfitted points. It is also a nonlinear method as opposed to other existing methods, which means that the relationship between the guiding points are considered nonlinearly unlike typical decomposition methods that obtain piece-wise linear functions. Further research for improving the accuracy of the degree of smoothness should be completed in practical applications.

## 2. Domain Decomposition and the Mesh-free Methods

The domain decomposition method is related to surface reconstruction. The reason is that we need to first determine the region that a sample point affects. The most natural decompositions include the Voronoi diagram and its dual correspondence, the Delaunay triangulation. Based on a set of sample points in a region, called sites, the Voronoi diagram partitions the domain into sub-regions where each point of the domain is categorized into a subset containing a site. The point is then considered closer to *this* site than to any other site. Its Delaunay triangulation is used to link the sites that have adjacent cells.

No other way can be better than this mechanism in general (if it is a must for a decomposition). The sub region that is a convex maintains a number of good properties. Many real world examples are available including the popular nearest neighbor method. Mount gives an excellent explanation of these questions in his lecture notes [22].

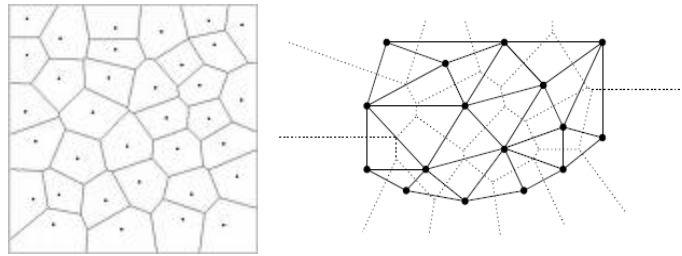

Fig 2. (a) A Voronoi diagram; (b) Relationship of a Voronoi diagram and its Delaunay triangulation [22].

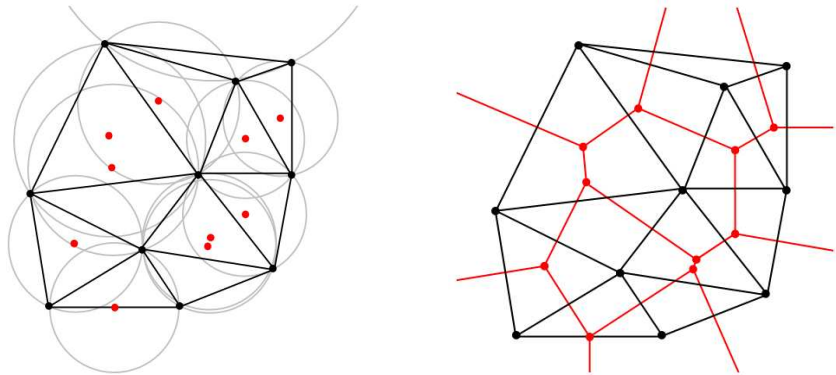

Fig 3. Delaunay triangulation vs. Voronoi diagram from Wikipedia
http://en.wikipedia.org/wiki/Delaunay_triangulation

Another way to define a Delaunay triangulation is to do a triangular decomposition on the sites such that no site is inside the circumcircle of any triangle. "The Delaunay triangulation with all the circumcircles and their centers (in red). Connecting the centers of the circumcircles produces the Voronoi diagram (in red)." (http://en.wikipedia.org/wiki/Delaunay_triangulation.)

A current popular data fitting method other than domain decomposition is called the moving least squares method (MLS) [17][3]. This method depends on the local region to fit the linear function locally and then applies the fitting to all "locations." The problem with this method is that it uses a weight function to determine the contribution from each "site" (sample points) to the unknown points. The function is usually dependent on different cases, so it is not fully automatic. Still, the local sample points determined for the moving point are the key. Using all or part of the region will cause the weight to be dependent on the individual cases. Even though, MLS does not need domain decomposition, it is very similar to this method in terms of local point determination.

Historically, MLS was studied in the 1970s [24][21] and used in graphics by Liven for point-set surfaces [18][19]. Amenta and Kil then used the disk to approximate the methods [2].

Recently, Fefferman proved the existence of a smooth function up to a certain smoothness degree order based on the refinements of grid points in Euclidean space applicable even in extended dimensions [13]. Fefferman et al described a method that is based on the local Whitney's jet problem and has some similarities with MLS.  Fefferman used the partition of unity for smoothing the "boundary" of each Whitney's jet.  The above methods are mainly linear methods or based on linear methods.

For pure refinement methods, with triangulation, if we know the function value at the vertex, we can insert three more middle points on each edge. Then we can use the mean of two end points to find the function value of the middle points. Therefore we have a continuous interpolation.   This is very similar to smooth color interpolation in graphics.

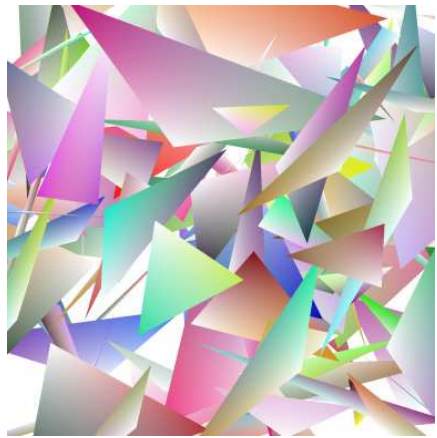

Fig 4. The interpolation for triangles by three corner points. http://actionsnippet.com/?p=2329

In fact, most uses of the practical methods to make a smooth surface on a manifold are in computer graphics. In 1978, Catmull and Clark discovered a recursive method to getting the smooth refinement of meshes [12]. The advantage of this method is that it allows "high resolution Mesh modeling without the need to save and maintain huge amounts of data." Therefore it is a refinement method.  It usually applies to ordinary shapes such as triangles, circles, or rectangles. Other methods for the same purpose use local rectangular or circular refinements [23][20]. However, these methods are not true interpolation and extension methods. They are designed for refinements meaning that we have to have relatively good existing decompositions or meshes before applying these methods.

Some algorithms are good for computer graphics but may not be practical for general data processing.  If we consider an example where the sample data points are located on one or two lines, we would not be able to fit the data appropriately. It is important to note that this case is very common in the real world.

### 3. Criticizing the Traditional Methods

In general, the classical geometric method uses domain decomposition to assist in obtaining the numerical interpolation of a domain that only contains random sample points.  This is because we need to determine the region that a sample point affects. Therefore, the Voronoi decomposition and its dual correspondence, the Delaunay triangulation, are the most popular in practice [26][22][1]. The Voronoi decomposition provides a piece-wise linear approximation to a manifold. This is also regarded as the classical discrete (data) reconstruction method and can be used with piecewise polynomials.

For example, when we have a Delaunay triangulation based on guiding points, we could define a Spline function since we know all the values of the vertices or we could add a PDE constraint on the sample points using the finite element method. As a result, some researchers would believe we have already found the perfect solution to the problem.

However, this is not the case in the real world. When we use the Voronoi decomposition for sample points, we are simultaneously assuming that the guiding points are "linearly separable." However, if we know that the unfitted function is not linear, we can conclude that the guiding points cannot be easily separated linearly or that we do not know the separation contour curve based on the guiding points if they are not dense enough. Euclidean space is a Jordan space because it is dense. If the guiding points are not dense and are randomly arranged, then the Voronoi decomposition or any other triangulation will generate false domain decompositions. This is the situation in Fig 5.

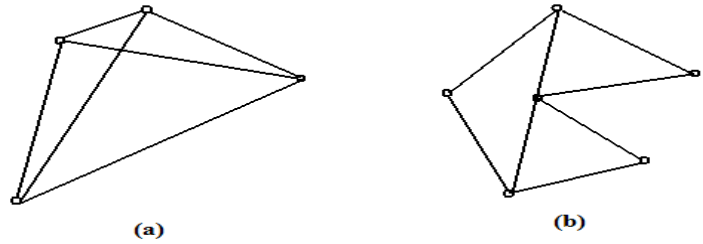

Fig 5. Valid cases but classical methods cannot be applied. (a) A nonlinear (non-Jordan) case. The point connections are not linearly separable. (b) Non-regular or improper mesh case; the Finite Element Method does not allow this type of decomposition.

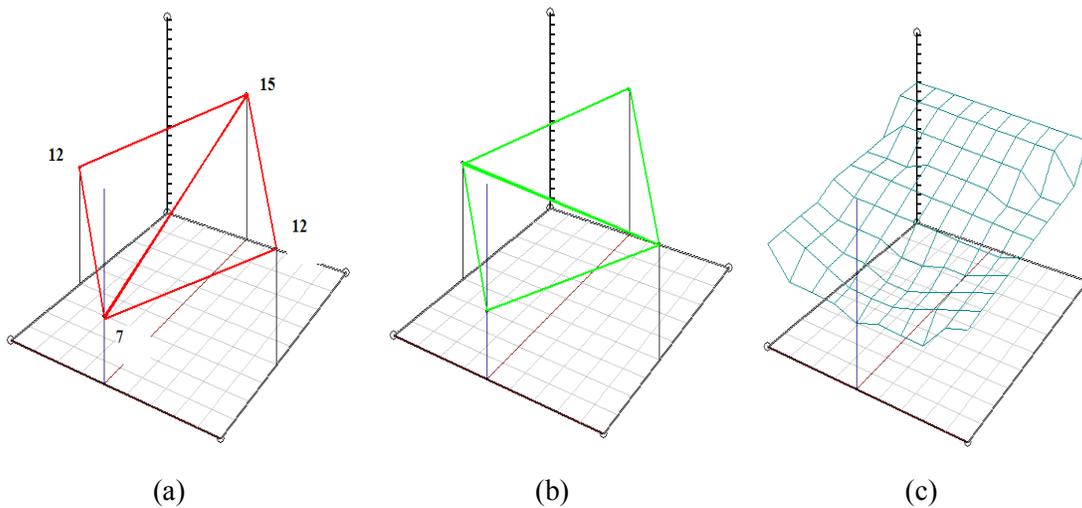

(a) (b) (c)

Fig 6. (a) (b) Two piecewise linear interpolations but we do not know which one is correct. (c) The GVF interpolation result shows a reasonable non-linear fitting.

If we have four sample points in the domain, then we will have two different piece-wise linear interpolations. The example is shown in Fig 6. One could say that averaging the two functions in Fig 6(a) and (b) would result in a better function. However, if we have five sample points, then we would have 10 different piece-wise linear interpolations (See Fig 7.). For six sample points, we may have more than 30 piece-wise linear functions. Taking the average value for each fitted point would be reasonable. However, the time complexity would be exponential.

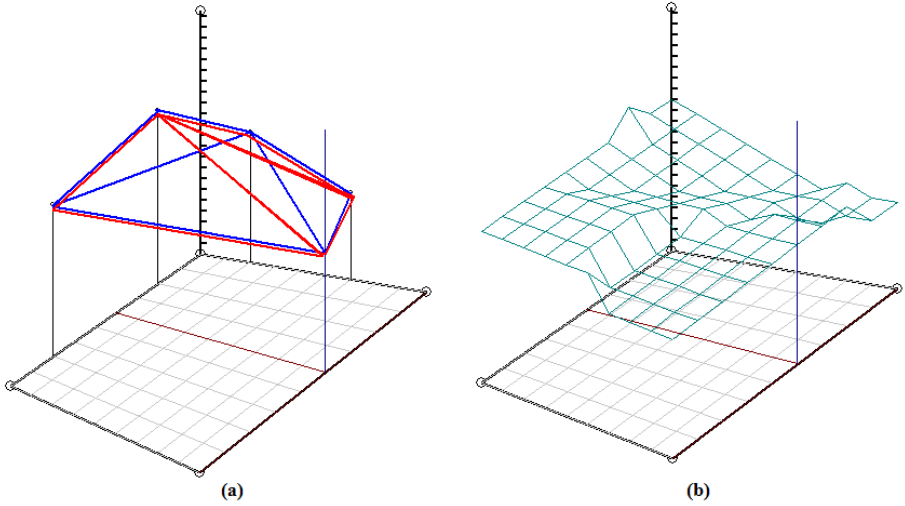

Fig 7. A case with five sample points will result in 10 piece-wise linear interpolations: (a) Two interpolations, (b) A GVF interpolation.

To avoid having too many possible cases in domain decomposition methods, mesh-free methods are proposed [3]. Moving least squares (MLS) is one of the most popular mesh-free methods in engineering and computer graphics. General speaking, MLS is good for many real world problems. The philosophy behind MLS is also one of the most efficient; it uses a tangent plane (usually a disk) to move to all unknown points with minimal errors. The drawback is determining the size of the disk. Depending on the situation, we may need to dynamically change its size. In order to consider more sample points, a weighted least square measurement is biased towards the region around the point at which the reconstructed value is requested [3][17][2][19]. Wikipedia provided a clear definition as follows:

"Consider a function $f : \mathbb{R}^n \rightarrow \mathbb{R}$ and a set of sample points $S = \{(x_i, f_i) \mid f(x_i) = f_i\}$ where $x_i \in \mathbb{R}^n$ and the $f_i$'s are a real numbers. Then, the moving least square approximation of degree $m$ at the point $x$ is $\tilde{p}(x)$ where $\tilde{p}$ minimizes the weighted least-square error

$$\sum_{i \in I}(p(x) - f_i)^2 \theta(\|x - x_i\|)$$

over all polynomials $p$ of degree $m$ in $\mathbb{R}^n$. $\theta(s)$ is the weight and it tends to zero as $s \rightarrow \infty$. In the example $\theta(s) = e^{-s^2}$." (http://en.wikipedia.org/wiki/Moving_least_squares)

An easy example could make MLS fail. Consider the Dirichlet problem: given a boundary of sample points, you are asked to fill the 3D surface. How would the weight function be chosen? MLS is a reasonable and practical method, but it requires a dense sample data set. Even though a Gaussian distribution weighted function was selected in general cases. If the surrounding region considered is too small, we would not have enough guiding points. This is why sometimes in artificial intelligence, k-nearest neighbors are required, but we would then need to determine the value of k. In addition, because there is no mesh, for a thin manifold, MLS may fail to determine the neighbor points by using Euclidean distance.

That is to say, for any fitting problem on a manifold, a system that recodes the adjacencies must be defined before the interpolation or fitting. Mesh-free is possible, however, adjacency relation free is not possible for such a calculation. That is why we can say that MLS is not a general-purpose interpolation method.

MLS is variation of Voronoi map, which is based on a measurement of a predefined weight function that can consider more guiding points in the calculation based on a weight that is decreasing by the distance to the point we are considering. It is a good idea for some applications, but would not work in all cases. The principle is still based on the linearity of domain separability. Even though there is a nonlinear weight function added, philosophically, the method is still "linear" since the importance decreases by a function of distance. In fact, MLS cannot automatically choose the weight functions.

Both of these methods only apply to cases with random but evenly distributed sample points.

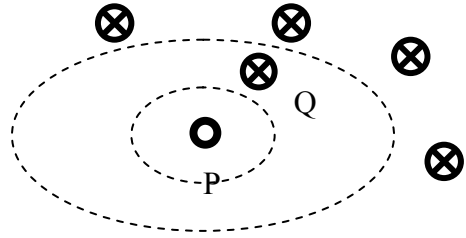

Fig 8. P is the moving point. Dotted ovals are the highly weighted areas in which the sites are contributing greatly to the value of point P. The cross X marks the sites.

The value of P might only be dependent on Q in Fig 8. However the outer four points will contribute evenly to P. In such a case, it would be easy to set up a sample point set to generate a false interpolation result.

For some P there are only one or two sites in the closer neighborhood. How do we determine the case? We may need the artificial intelligence method to dynamically determine the theta function. We can see that MLS does not depend on mesh but on circles so there must be dense samples on manifolds. Therefore, the existing practical methods for "linear" decomposition or only smooth refinements will not be accurate. A true interpolation method is needed for random data sample points.

## 4. Gradually Varied Functions

For a manifold, there are two types of meshes: 1) Domain meshes, and 2) Guiding point's meshes. A manifold can be defined by CW-complexes in topology. We could assume that a fine decomposition of simplexes can represent the manifold, i.e. the linear decomposition (piece-wise linear representation).

The guiding point mesh, if there is a good one, is not as dense as the domain mesh. Since it may only contain a few points, we cannot assume that there is linearity among the guiding points, regardless of whether we could assume linearity among domain points. If we have to use the triangulations in guiding points, there is a reasonable way to do this. We could find all the possible simplex decompositions, and then find the average to approximate the original. However, it would cost exponential time to accomplish this task by thinking about the number of different forms of simplex decompositions for n points.

The way that may overcome or solve the problem is to find a nonlinear reconstruction method. The new method would not use domain decomposition. The gradually varied function was proposed in 1989 to solve a problem with the random data sample points [4][5][6][7] and fortunately is a non-linear interpolation method.

Gradual variation is a discrete method that can be built on any graph and the gradually varied surface is a special discrete surface. We now introduce this concept.

Concept of Gradual Variation: Let function f: D→{A1, A2,…,An}. Given that a and b are adjacent in D implies f(a)=f(b), or f(b) =A(i-1) orA(i+1) when f(a)=Ai. Point (a,f(a)) and (b,f(b)) are then said to be gradually varied. A 2D function (surface) is said to be gradually varied if every adjacent pair is gradually varied.

Discrete Surface Fitting: Given J⊆D, and f: J→{A1,A2,…An} decide if there is an F: D→{A1,A2,…,An} such that F is gradually varied where f(x)=F(x), for every x in J.

**Theorem** (Chen, 1989) [4] The necessary and sufficient conditions for the existence of a gradually varied extension F is: for all x, y in J, $d(x,y) \geq |i-j|$, f(x)=Ai and f(y)=Aj, where d is the distance between x and y in D.

Assume a domain is linearly separable, such as the Jordan domain. If we need a local region, we can first fit the boundary of the region and then fit the inner points. Our discrete algorithm does not limit the order of the reconstruction. So the gradually varied function fitting is both a local and global method [4][5][6]. In the above sections, we show several examples using GVF and we have conducted more real data processing in [10].

## 5. Smooth Digital-Discrete Surface Reconstruction

Can gradually varied functions be used to obtain smooth functions? Chen and Adjei proposed a method to solving this problem in 2004 [11]. Fefferman used Whitney's Jets (a refinement method) and linear programming to ensure the existence of the $C^{(m)}$ extension [13]. In fact, for the continuous function ($C^{(0)}$) case, some of these treatments are believed equivalent or similar to gradually varied functions under the Lipschitz condition, the most important case in many real world problems. Recently, Shvartsman discussed the extension in Sobolev space using Lipschitz conditions [15], a method for $C^{(1)}$ refinements was discussed in [12], and Fefferman still uses the local refinement method in his research.

In [8], we directly deal with smooth functions on the (same) piecewise linear manifold or non-Jordon graph. This paper presents applications of using this method in many different cases including rectangle domains and closed surfaces. We also have applied this method to groundwater flow equations [9]. The new algorithm is based on the gradually varied functions and finite difference methods. To find a theoretical extension, GVS can also be a powerful method [8]. Fig 9 shows the result of smooth fittings.

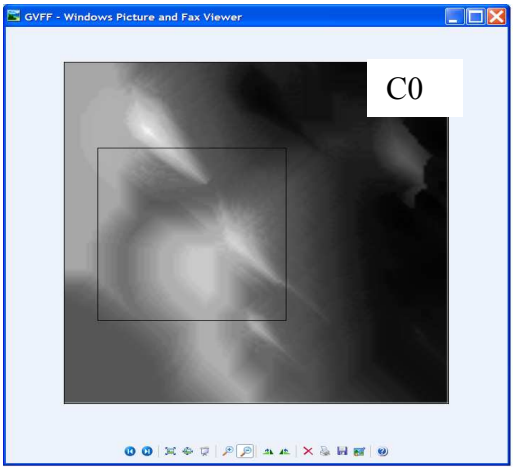 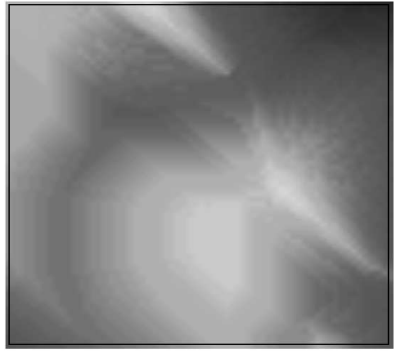

(a) $C^0$

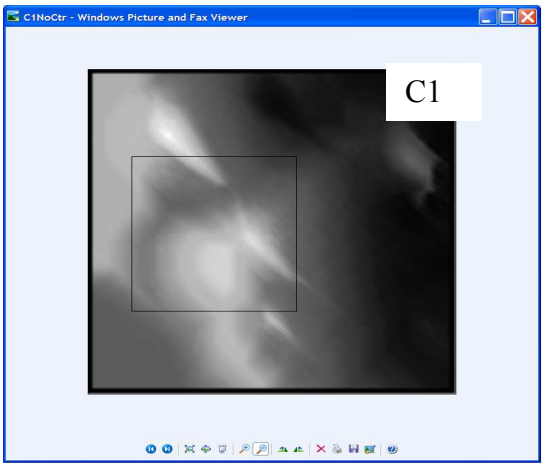 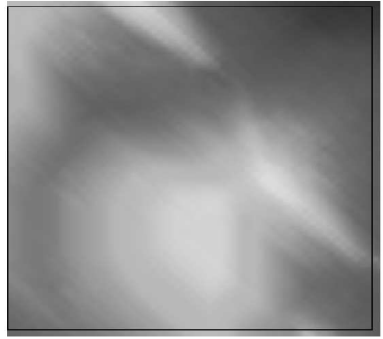

(b) $C^1$

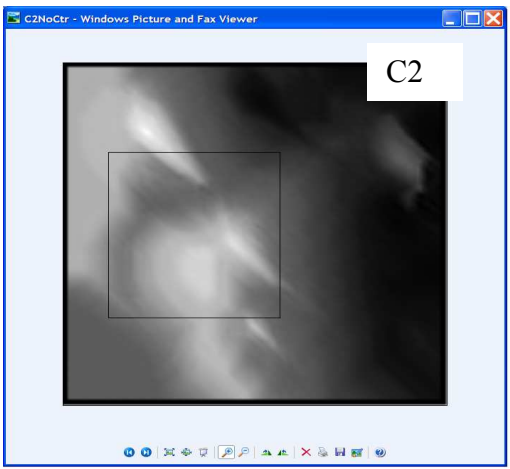 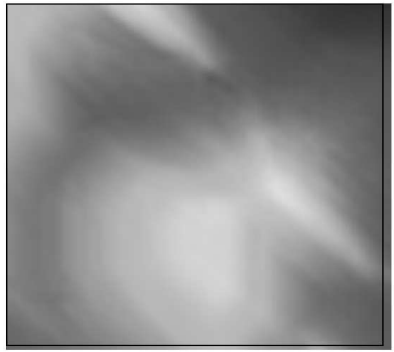

(c) $C^2$

Fig 9. The result of interpolation based on 29 sample points. (a) The first is a "continuous" surface; (b) The result has a continuous "first derivative."  (c) The result has a continuous "second derivative."

For the function on manifolds, we have the following results [10]. The guiding points form a cycle, and we can see that the GVF fitting and harmonic fitting are similar. Since the harmonic function in this case will converge after a certain number of iterations, it validates the result from GVF.

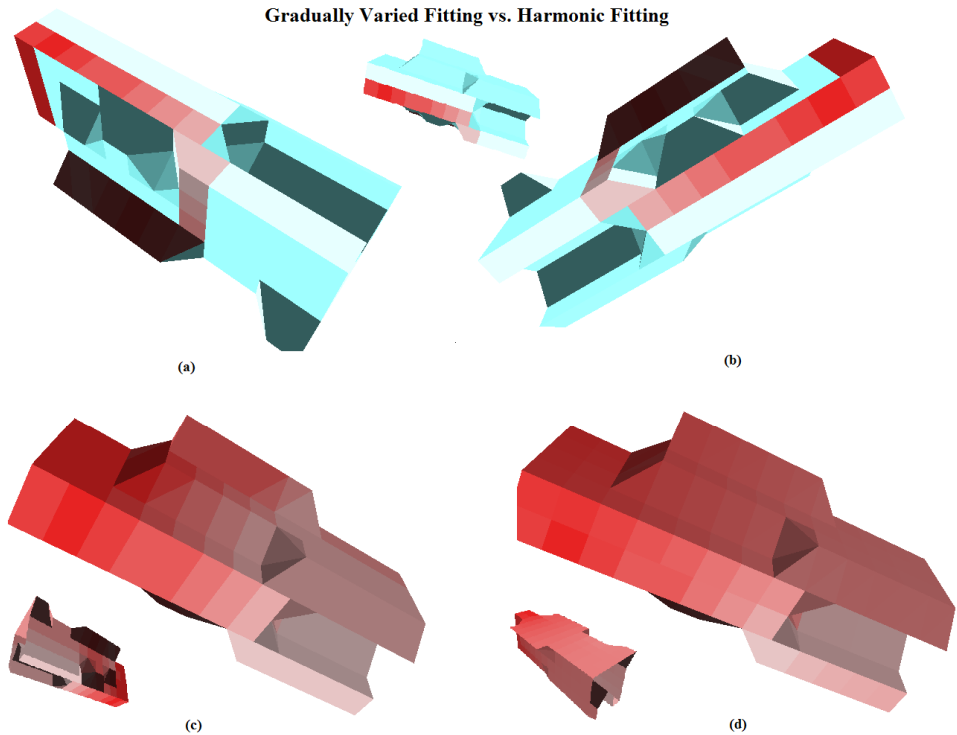

Gradually Varied Fitting vs. Harmonic Fitting: (a) The selected cells form a boundary curve that is gradually varied. (b) Another view of the guiding points. (c) The gradually varied fitting (GVF) result. (d) The Harmonic fitting based on GVF (100 iterations).

## Appendix: Remarks

For continuity in discrete space, we have discrete harmonic functions, which has been known for some time by mathematicians. For digital space (mostly for digital images), in 1985, Chen proposed lambda-connectedness. In 1986, A. Rosenfeld presented discretely continuous functions. In 1989, Chen simplified lambda-connectedness to gradual variation. Also in 1989, J.M. Steele used the Lipschtz condition to check for image "smoothness," which was in fact the same as the continuity in Chen's paper. However, Steele did not mention the reconstruction problem, which should be the more important part to him as a mathematician at Princeton. He moved on to other problems, never writing a second paper in this same category, as he told the author. As an expert also in computational geometry, Steele could have helped the establishment a great deal if he had continued.

On the other hand, Rosenfeld, who was regarded as the father of computer vision, realized the reconstruction problem for images. However for image processing, the reconstruction problem was not a big concern compared to the segmentation problem. In general, Rosenfeld was a grandmaster for generating new concepts. Perhaps he was not very interested in detailed results or he had met some difficulties; he seemed to have chosen not to continue pursuing this problem. Rosenfeld had many students and young collaborators. None of them seriously picked up this direction for further research. To be frank, his related statement and proof on reconstruction were not correctly stated in his 1986 paper. Interestingly enough, Chen began to learn digital topology after studying Rosenfeld's earlier paper on digital surfaces. Through the research papers, Rosenfeld was a great advisor in inspiring many young scientists, Chen being one of them. Chen was well trained in computational complexity and discrete mathematics, including graph theory. When he worked in geophysical prospecting for petroleum, due to the placement of irregular sampling lines, Chen realized that the reconstruction of surfaces was just as important as the segmentations for data processing. Since he had knowledge of graph theory at the time in 1988-1989, he could come up with a very elegant proof for the main theorem.

It is obvious that the problem of fitting random finite points to a surface has great importance in both theory and practice. Recently, Annals of Mathematics regarded as the best mathematical journal in the world published three papers by Fields Medalist C. Fefferman. This is another indication of the importance of the problem.